# Optimal control of impulsive Volterra equations with variable impulse times

S. A. Belbas  
Mathematics Department  
University of Alabama  
Tuscaloosa, AL 35487-0350.  
USA.

W. H. Schmidt  
Institut für Mathematik und Informatik  
Universität Greifswald  
D-17487 Greifswald.  
Germany.

e-mail: sbelbas@bama.ua.edu

e-mail: wschmidt@uni-greifswald.de

Abstract. We obtain necessary conditions of optimality for impulsive Volterra integral equations with switching and impulsive controls, with variable impulse time-instants. The present work continues and complements our previous work on impulsive Volterra control with fixed impulse times.

1. Introduction.

The traditional theory of optimal control has been developed mostly for differential equations (ordinary or partial), delay-differential equations, and differential inclusions.

An equally important class of controlled systems consists of controlled Volterra integral equations. In fact, the class of controlled Volterra integral equations also encompasses, as particular cases, all controlled ordinary differential equations, as well as controlled integrodifferential equations. By comparison to most other classes of control systems, relatively little is known about controlled Volterra integral equations. It is part of our intentions to attract more interest to this important category of optimal control problems. Volterra equations describe systems with continuously distributed memory over the entire past of the system; consequently, they have features that are substantially different from those of memoryless systems (i.e. ordinary or partial differential equations and differential inclusions), and also very different from those of systems with concentrated memory effects (i.e. delay differential equations, with either constant or variable delays).

Another important direction in optimal control theory concerns the class of impulsive systems, with either continuous, or impulsive, or switching controls, or any combination of these three types of controls. Here, we must distinguish between the concepts of impulsive systems, which concerns the form of state dynamics, and the concept of switching controls, which concerns the class of admissible control actions. For example, in the classical theory of time-optimal control, the so-called "bang-bang" controls are switching controls, but the state trajectory has no impulses. A system may involve controls that simultaneously induce switchings in the state dynamics as well as impulses in the state trajectory. The classification of state dynamics and control actions is also reflected in the cost functional, which can include continuous, switching, and impulsive terms.

Many systems involve memory effects and can be modelled by Volterra integral equations. For example, the spread of epidemics with general (not necessarily exponential) distribution of infection times is modelled by Volterra integral equations; the evolution of a set of competing biological species can be modelled by a set of ODEs of Volterra-Lotka type, but when memory effects are included, then the model becomes a system of integro-differential equations, which can be reduced to a system of Volterra integral equations; in Economics, the evolution of capital stock under an investment strategy can be described by a Volterra integral equation.
We are interested in controlled Volterra integral equations with impulsive actions of the control. Impulsive controls arise in a variety of applications, including discrete decisions in finance problems, harvesting of biological populations at discrete time-instants, pulse vaccination in problems of control of epidemics, etc. The theory of optimal control of Volterra integral equations with continuous controls has been developed in [S1] and certain other papers that are referenced in [BS]. Related control problems for differential equations have been studied in [B]. By contrast, relatively little is known about impulsive controlled Volterra integral equations. We mention that the paper [BS] contains results on optimal control of impulsive Volterra integral equations with fixed impulse times. Control of Volterra equations with switching controls and variable switching times, but without impulses in the state dynamics, has been studied in [S2, S3, S4]. The general problem of impulsive optimal control of Volterra integral equations involves two kinds of control decisions: (i) the optimal time-instants at which impulses are applied, and (ii) the optimal values of the impulsive control actions at the optimal time instants. Of course, these two aspects are not independent from each other, but rather the two types of optimal decisions are coupled into one optimal control problem.

In this paper we derive necessary conditions arising out of variations of the optimal impulse time-instants. The results presented herein complement those of [BS]. Thus, with the results of the present paper in addition to those of [BS], we have achieved a complete set of necessary conditions for the general problem of optimal control of impulsive Volterra integral equations with variable impulse times.

## 2. Statement of the problem.

We consider an impulsive controlled Volterra integral equation with piecewise constant controls. The time-horizon is [0, T]. The impulsive times are

$$0 < \tau_1 < \tau_2 < ... < \tau_N < T.$$

We set $\tau_0 := 0$, $\tau_{N+1} := T$, and $S_N := \{(\tau_1, \tau_2, \cdots, \tau_N): 0 < \tau_1 < \tau_2 < ... < \tau_N < T\}$.
The control u(t) is defined as

$$u(t) = a_i \text{ for } \tau_i < t < \tau_{i+1}, \ i = 0,1,2,..., N.$$

Each of the constants $a_i$ takes values in a closed bounded nonempty set in $R^m$ that is equal to the closure of its own interior.

Let $y(\cdot):[0,T] \mapsto R^n$ be the state function.

We set $\boldsymbol{\tau} := \{\tau_1, \tau_2,..., \tau_N\}$, $\mathbf{a} := \{a_0, a_1,..., a_N\}$, $\mathbf{y} := \{y(\tau_1^-), y(\tau_2^-),..., y(\tau_{N+1}^-)\}$. For each t in [0, T], we set

$$\tau_{(t)} := \{\tau_i : i \geq 1, \tau_i < t\}, \ y_{(t)} := \{y(\tau_i^-): i \geq 1, \tau_i < t\}, \ a_{(t)} := \{a_i : i \geq 0, \ \tau_i < t\}.$$

The symbol $y(\tau_i^-)$ denotes the one-sided limit $\lim_{t \to \tau_i^-} y(t)$ ; analogously,

$$y(\tau_i^+) := \lim_{t \to \tau_i^+} y(t).$$

The state dynamics is described by the impulsive Volterra equation

$$y(t) = y_0(t) + \int_0^t f(t, s, y(s), u(s))\,ds + g(t, \tau_{(t)}, y_{(t)}, a_{(t)}) \tag{2.1}$$

The cost functional to be minimized is

$$J := \int_0^T F(t, y(t), u(t))\,dt + G(\boldsymbol{\tau}, \mathbf{y}, \mathbf{a}) \tag{2.2}$$

We note that the state dynamics can be written as

$$y(t) = y_0(t) + \sum_{i: i \geq 0, \, \tau_i < t} \int_{\tau_i}^{t \wedge \tau_{i+1}} f(t, s, y(s), a_i) ds + g(t, \tau_{(t)}, y_{(t)}, a_{(t)}) \qquad (2.3)$$

(where $a \wedge b := \min(a, b)$) and the cost functional J can also be written in the form

$$J = \sum_{i=0}^{N} \int_{\tau_i}^{\tau_{i+1}} F(t, y(t), a_i) dt + G(\boldsymbol{\tau}, \mathbf{y}, \mathbf{a}) \qquad (2.4)$$

Because of the causal structure of Volterra equations, the value y(t) of the solution y at time t depends on $t, \tau_{(t)}, a_{(t)}$.

Eq. (2.3), in addition to being useful in other respects, makes explicit the impulsive nature of the state dynamics, since it plainly implies that

$$y(\tau_i^+) - y(\tau_i^-) = g(\tau_i, \{\tau_j : 1 \leq j \leq i\}, \{y(\tau_j^-) : 1 \leq j \leq i\}, \{a_j : 0 \leq j \leq i\}) -$$
$$- g(\tau_i, \{\tau_j : 1 \leq j \leq i-1\}, \{y(\tau_j^-) : 1 \leq j \leq i-1\}, \{a_j : 0 \leq j \leq i-1\}) \quad .$$

We postulate the following conditions on the functions that appear in (2.1) and (2.2):

(1). $y_0(\cdot)$ is continuous on [0, T].
(2). f is continuous with respect to all its arguments, and continuously differentiable with respect to y and u, with bounded derivative with respect to y.
(3). g is continuous in all its arguments, and continuously differentiable with respect to $\tau_{(t)}, y_{(t)}, a_{(t)}$ for all t in [0, T].
(4). The function F is continuous with respect to all its arguments, and continuously differentiable in y and u.
(5). The function G is continuously differentiable in $\boldsymbol{\tau}, \mathbf{y}, \mathbf{a}$.

Existence and uniqueness of solutions of the impulsive Volterra equation, and constructive methods of approximating the solutions, follow from [BS].
A note on notation: partial differentiation with vector arguments denotes a Jacobian matrix of appropriate dimensions, for example $\frac{\partial f}{\partial y}$ or $f_y$ denotes the matrix with the entry in the $\alpha$-th row and $\beta$-th column equal to $\frac{\partial f_\alpha}{\partial y_\beta}$ where $f = \text{col}[f_\alpha : 1 \leq \alpha \leq n]$, $y = \text{col}[y_\beta : 1 \leq \beta \leq n]$, and "col" denotes a column vector with the indicated components. In this notation, the gradient of a scalar-valued function is a row vector. Various quantities that appear below will be vectors (column or row vectors) or matrices of appropriate dimensions; the dimensions will not always be explicitly specified, but rather will be understood from the context, i.e. from the requirement of compatibility of dimensions in the sense of making the appropriate algebraic operations meaningful.

We note that the differentiability conditions are to be utilized in the variational analysis in the remaining of this paper; the existence and uniqueness of solutions of (2.1) can be proved under less stringent conditions.

## 3. Background on impulsive linear Volterra integral equations.

The variational equations arising out of impulsive nonlinear Volterra integral equations are linear (as equations with the variations as unknowns), and they are still impulsive Volterra equations. Consequently, in order to carry out the variational analysis for the problem formulated in section 2, we need results on linear impulsive Volterra equations, and in particular the suitable concept of the adjoint equation for a linear impulsive Volterra equation. We shall utilize certain results from [BS].

We consider the impulsive Volterra equation

$$y(t) = \eta(t) + \int_0^t K(t,s)y(s)\,ds + \sum_{i:0<\tau_i<t} \lambda_i(t)y(\tau_i^-)$$

--- (3.1)

The functions $\eta(.), \lambda_i(.)$ are continuous on each interval $(\tau_j, \tau_{j+1})$ with finite limits as $t \to \tau_j^{\pm}$ for each $\tau_j$, and $K(.,.)$ is continuous in (t, s) for $0 \le s \le t, t \in (\tau_j, \tau_{j+1})$ with finite limits as $t \to \tau_j^{\pm}$ for each $\tau_j$.

We define the $N \times N$ array $\mathbf{\Lambda} \equiv [\Lambda_{ij}]$ by $\Lambda_{ij} := \lambda_j(\tau_i^-)$, $j = 1, 2, .., i-1$, $i = 2, 3, ..., N$, and $\Lambda_{ij} = 0$ for all other values of i, j. Each $\Lambda_{ij}$ is an $n \times n$ matrix. Then it follows from (3.1) that

$$y(\tau_i^-) = \eta(\tau_i^-) + \int_0^{\tau_i} K(\tau_i^-, s)y(s)\,ds + \sum_{j=1}^{N} \Lambda_{ij} y(\tau_j^-)$$

--- (3.2)

The array $\mathbf{\Lambda}$ is lower-block-triangular with zero diagonal, thus $\mathbf{\Lambda}^N = \mathbf{0}$ (the product $\mathbf{\Lambda M}$ of two quantities $\mathbf{\Lambda}$ and $\mathbf{M}$ having the dimensions of $\mathbf{\Lambda}$ is defined as $\mathbf{\Lambda M} = \sum_{k=1}^{N} \Lambda_{ik} M_{kj}$), and the solution of (3.2) by the method of simple iterations gives

$$y(\tau_i^-) = \eta(\tau_i^-) + \int_0^{\tau_i} K(\tau_i^-, s)y(s)\,ds +$$
$$+ \sum_{j=1}^{i-1} \left(\mathbf{\Lambda} + \mathbf{\Lambda}^2 + \cdots + \mathbf{\Lambda}^{N-1}\right)_{ij} \left(\eta(\tau_j^-) + \int_0^{\tau_j} K(\tau_j^-, s)y(s)\,ds\right)$$

--- (3.3)

We set

$$\Gamma_{ij} := \left(\Lambda + \Lambda^2 + \cdots + \Lambda^{N-1}\right)_{ij}$$

--- (3.4)

We remark that each $\Gamma_{ij}$ can be expressed in terms of increasing paths from j to i. An increasing path from j to i is a collection of indices $k_1, k_2, \ldots, k_\alpha$, out of the set of indices that correspond to impulsive time-instants, with the property $j < k_1 < k_2 < \cdots < k_\alpha < i$. Then $\Gamma_{ij} = \sum \Lambda_{ik_\alpha} \Lambda_{k_\alpha k_{\alpha-1}} \cdots \Lambda_{k_1 j}$, where the summation is taken over all increasing paths from j to i.

By using the notation defined in (3.4), we can write (3.1) in the form

$$y(t) = \eta(t) + \sum_{i:0<\tau_i<t} \sum_{j=1}^{i-1} \lambda_i(t) \Gamma_{ij} \left( \eta(\tau_j^-) + \int_0^{\tau_j} K(\tau_j^-, s) y(s) \, ds \right) + \int_0^t K(t,s) y(s) \, ds$$

--- (3.5)

Thus we have a new Volterra integral equation of the form

$$y(t) = \tilde{\eta}(t) + \int_0^t \tilde{K}(t,s) y(s) \, ds;$$

$$\tilde{\eta}(t) := \eta(t) + \sum_{i:0<\tau_i<t} \sum_{j=1}^{i-1} \lambda_i(t) \Gamma_{ij} \eta(\tau_j^-), \quad \tilde{K}(t,s) = K(t,s) + \sum_{i:0<\tau_i<t} \sum_{j=1}^{i-1} \lambda_i(t) \Gamma_{ij} K(\tau_j^-, s)$$

--- (3.6)

The solution of (3.6) can be expressed in terms of a resolvent kernel $R(t,s)$ defined by

$$R(t,s) = \sum_{n=1}^{\infty} \tilde{K}_n(t,s); \quad \tilde{K}_1(t,s) := \tilde{K}(t,s), \quad \tilde{K}_{n+1}(t,s) := \int_s^t \tilde{K}_n(t,\sigma) \tilde{K}(\sigma,s) \, d\sigma, \quad n \geq 1$$

--- (3.7)

and y(.) is given by

$$y(t) = \tilde{\eta}(t) + \int_0^t R(t,s)\tilde{\eta}(s)\,ds$$

--- (3.8)

It follows, as in [BS], that the function

$$z(t) = \tilde{\eta}(t) + \int_t^T \tilde{\eta}(s)R(s,t)\,ds$$

--- (3.9)

satisfies

$$z(t) = \tilde{\eta}(t) + \int_t^T z(s)K(s,t)\,ds + \int_t^T \sum_{i:0<\tau_i<s} \sum_{j=1}^{i-1} z(s)\lambda_i(s)\Gamma_{ij}K(\tau_j^-,t)\,ds$$

--- (3.10)

By changing the order of summations and integration, and keeping in mind that $\tau_i > \tau_j$ for the summations shown in (3.10) and that $K(\tau_j^-, t) = 0$ for $t \geq \tau_j$ and $\lambda_i(s) = 0$ for $s \leq \tau_i$, we can write (3.10) in the equivalent form

$$z(t) = \tilde{\eta}(t) + \int_t^T z(s)K(s,t)\,ds + \sum_{j:t<\tau_j<T} \sum_{i=j+1}^N \left( \int_{\tau_i}^T z(s)\lambda_i(s)\,ds \right) \Gamma_{ij} K(\tau_j^-, t)$$

--- (3.11)

Eq. (3.11) is the adjoint equation to the linear impulsive Volterra integral equation (3.1).

### 4. Variational equations for the impulsive times.

This section concerns the variational equations of the state dynamics (3.1) with respect to variations in the impulse times. We calculate the following derivatives:

$$\partial_j y(t) := \frac{\partial y(t; \tau_{(t)}, a_{(t)})}{\partial \tau_j}, \quad \tau_j < t, \ t \notin \tau; \ 1 \le j \le N;$$

$$\partial_j y(\tau_\ell^-) := \frac{\partial y(\tau_\ell^-; \tau_{(\tau_\ell)}, a_{(\tau_\ell)})}{\partial \tau_j}, \quad j < \ell \le N+1;$$

$$D_j y(\tau_j^-) := \frac{\partial y(\tau_j^-; \tau_{(\tau_j)}, a_{(\tau_j)})}{\partial \tau_j}, \quad 1 \le j \le N$$

--- (4.1)

By differentiating (2.3), we obtain, for $t \notin \tau$ with $t > \tau_j$,

$$\partial_j y(t) = f(t, \tau_j, y(\tau_j^-), a_{j-1}) - f(t, \tau_j, y(\tau_j^+), a_j) +$$

$$+ \sum_{i:\ \tau_j \le \tau_i < t} \int_{\tau_i}^{t \wedge \tau_{i+1}} \frac{\partial f(t, s, y(s), a_i)}{\partial y} \partial_j y(s) \, ds +$$

$$+ \frac{\partial g(t,...)}{\partial \tau_j} + \frac{\partial g(t,...)}{\partial y_j} D_j y(\tau_j^-) + \sum_{\ell:\ \tau_j < \tau_\ell < t} \frac{\partial g(t,...)}{\partial y_\ell} \partial_j y(\tau_\ell^-)$$

--- (4.2)

In order to calculate $D_j y(\tau_j^-)$, we use

$$y(\tau_j^-) = y_0(\tau_j) + \sum_{i=0}^{j-1} \int_{\tau_i}^{\tau_{i+1}} f(\tau_j, s, y(s), a_i) \, ds + g(\tau_j; \tau_{(\tau_j)}, y_{(\tau_j)}, a_{(\tau_j)})$$

--- (4.3)

from which, by differentiating with respect to $\tau_j$, we find

$$D_j y(\tau_j^-) = \frac{dy_0(\tau_j)}{d\tau_j} + f(\tau_j, \tau_j^-, y(\tau_j^-), a_{j-1}) + \int_0^{\tau_j} \frac{\partial f(\tau_j, s, ...)}{\partial \tau_j} ds + \frac{\partial g(\tau_j, ...)}{\partial \tau_j}$$

--- (4.4)

We remark that, according to our definition in section 2, the terms $\boldsymbol{\tau}_{(\tau_j)}, \mathbf{y}_{(\tau_j)}, \mathbf{a}_{(\tau_j)}$ that appear on the right-hand side of (4.3) (and thus also on the right-hand side of (4.4)) are given by

$$\boldsymbol{\tau}_{(\tau_j)} = \{\tau_1, \tau_2, ..., \tau_{j-1}\}, \ \mathbf{y}_{(\tau_j)} = \{y(\tau_1^-), y(\tau_2^-), ..., y(\tau_{j-1}^-)\}, \ \mathbf{a}_{(\tau_j)} = \{a_0, a_1, ..., a_{j-1}\}.$$

The derivatives $\partial_j y(\tau_\ell^-)$ for $\ell > j$ can be calculated in the same way as the derivatives $\partial_j y(t)$ with $t > \tau_j$:

$$\partial_j y(\tau_\ell^-) = f(\tau_\ell, \tau_j, y(\tau_j^-), a_{j-1}) - f(\tau_\ell, \tau_j, y(\tau_j^+), a_j) + \int_{\tau_j}^{\tau_\ell} \frac{\partial f(\tau_\ell, s, ...)}{\partial y} \partial_j y(s) ds +$$

$$+ \frac{\partial g(\tau_\ell, ...)}{\partial \tau_j} + \frac{\partial g(\tau_\ell, ...)}{\partial y_j} D_j y(\tau_j^-) + \sum_{\ell'=j+1}^{\ell-1} \frac{\partial g(\tau_\ell, ...)}{\partial \tau_{\ell'}} \partial_j y(\tau_{\ell'}^-)$$

--- (4.5)

Eq. (4.2) is of the type

$$\partial_j y(t) = \eta_j(t) + \int_{\tau_j}^t K(t, s) \partial_j y(s) ds + \sum_{\ell: \tau_j < \tau_\ell < t} \lambda_\ell(t) \partial_j y(\tau_\ell^-)$$

--- (4.6)

where

$$K(t, s) := \frac{\partial f(t, s, y(s), u(t))}{\partial y};$$

$$\eta_j(t) := f(t, \tau_j, y(\tau_j^-), a_{j-1}) - f(t, \tau_j, y(\tau_j^+), a_j) + \frac{\partial g(t; ...)}{\partial \tau_j} + \frac{\partial g(t; ...)}{\partial y_j} D_j y(\tau_j^-);$$

$$\lambda_\ell(t) := \frac{\partial g(t; ...)}{\partial y_\ell}$$

--- (4.7)

According to the results of section 3, $\partial_j y(t)$ can be represented as

$$\partial_j y(t) = \tilde{\eta}_j(t) + \int_{\tau_j}^{t} R_j(t,s)\tilde{\eta}_j(s)\,ds \text{ for } t > \tau_j;$$

$$\tilde{\eta}_j(t) = \eta_j(t) + \sum_{i:\tau_j < \tau_i < t} \sum_{\ell=j}^{i-1} \lambda_i(t)\Gamma_{i\ell}\,\eta_j(\tau_\ell^-)$$

--- (4.8)

where the resolvent kernel $R_j(t,s)$ is constructed as in section 2.

Our exposition will be facilitated by introducing the oscillation operator $\Omega$, defined as

$$\Omega\varphi(\tau_i, u(.)) := \varphi(\tau_i^-, a_{i-1}) - \varphi(\tau_i^+, a_i)$$

--- (4.9)

The function $\varphi$ in (4.9) may depend on additional variables, other than the ones shown in (4.9), but the operator $\Omega$ is applied with respect to the variables indicated in (4.9). Thus, for instance,

$$\Omega f(t, \tau_j; y, u(.)) := f(t, \tau_j, y(\tau_j^-), a_{j-1}) - f(t, \tau_j, y(\tau_j^+), a_j)$$

--- (4.10)

With this notation, we have

$$\eta_j(t) = \Omega f(t, \tau_j; y, u(\cdot)) + \frac{\partial g(t,...)}{\partial \tau_j} + \frac{\partial g(t,...)}{\partial y_j} D_j y(\tau_j^-)$$

--- (4.11)

## 5. The impulsive extremum principle.

We want to define a concept of co-state for the variation of the cost functional J under variations for the impulse time, and we want to derive equations analogous to the Hamiltonian equations of ordinary control theory and an extremum principle for the problem under consideration.

We start by calculating the derivatives of the cost functional J with respect to the impulse times. We have

$$\frac{\partial J}{\partial \tau_j} = F(\tau_j, y(\tau_j^-), a_{j-1}) - F(\tau_j, y(\tau_j^+), a_j) +$$

$$+ \int_{\tau_j}^{T} \frac{\partial F(t, y(t), u(t))}{\partial y} \partial_j y(t) \, dt + \frac{\partial G(\boldsymbol{\tau}, \mathbf{y}, \mathbf{a})}{\partial \tau_j} + \frac{\partial G(\boldsymbol{\tau}, \mathbf{y}, \mathbf{a})}{\partial y_j} D_j y(\tau_j^-) +$$

$$+ \sum_{\ell=j+1}^{N+1} \frac{\partial G(\boldsymbol{\tau}, \mathbf{y}, \mathbf{a})}{\partial y_\ell} \partial_j y(\tau_\ell^-) =$$

$$= \Omega F(\tau_j; y, u(.)) + \int_{\tau_j}^{T} \frac{\partial F(t, y(t), u(t))}{\partial y} \partial_j y(t) \, dt +$$

$$+ \frac{\partial G(\boldsymbol{\tau}, \mathbf{y}, \mathbf{a})}{\partial \tau_j} + \frac{\partial G(\boldsymbol{\tau}, \mathbf{y}, \mathbf{a})}{\partial y_j} D_j y(\tau_j^-) + \sum_{\ell=j+1}^{N+1} \frac{\partial G(\boldsymbol{\tau}, \mathbf{y}, \mathbf{a})}{\partial y_\ell} \partial_j y(\tau_\ell^-)$$

--- (5.1)

By inserting the representation (4.8) for $\partial_j y(t)$ into (5.1), we obtain

$$\frac{\partial J}{\partial \tau_j} = \Omega F(\tau_j; y, u) + \int_{\tau_j}^{T} [F_y(t, y(t), u(t)) + \int_{t}^{T} F_y(s, y(s), u(s)) R_j(s, t) \, ds] \tilde{\eta}_j(t) \, dt +$$

$$+ \frac{\partial G(\boldsymbol{\tau}, \mathbf{y}, \mathbf{a})}{\partial \tau_j} + \frac{\partial G(\boldsymbol{\tau}, \mathbf{y}, \mathbf{a})}{\partial y_j} D_j y(\tau_j^-) +$$

$$+ \sum_{l=j+1}^{N+1} \frac{\partial G(\boldsymbol{\tau}, \mathbf{y}, \mathbf{a})}{\partial y_\ell} [\tilde{\eta}_j(\tau_\ell^-) + \int_{\tau_j}^{\tau_l} R_j(\tau_\ell^-, t) \tilde{\eta}_j(t) \, dt]$$

--- (5.2)

We define the co-states $p_j(t)$, $t > \tau_j$, by

$$p_j(t) := F_y(t, y(t), u(t)) + \int_t^T F_y(s, y(s), u(s)) R_j(s, t) ds$$

--- (5.3)

and the collection of Hamiltonians $h_j(t, y, u, p_j(.)), t > \tau_j$ by

$$h_j(t, y, u, p_j(.)) := F(t, y, u) + \int_t^T p_j(s) f(s, t, y, u) ds$$

--- (5.4)

Then, in view of the results of section 3, the co-states satisfy

$$p_j(t) = \frac{\partial}{\partial y} h_j(t, y(t), u(t), p_j(\cdot)) +$$

$$+ \sum_{\ell: t < \tau_\ell < T} \sum_{i=\ell+1}^N \left( \int_{\tau_i}^T p_j(s) \lambda_\ell(s) ds \right) \Gamma_{i\ell} f_y(\tau_\ell, t, y(t), u(t))$$

--- (5.5)

Now, (5.2) and (5.3) imply that the derivatives of the cost functional J with respect to the impulse times can be expressed as

$$\frac{\partial J}{\partial \tau_j} = \Omega F(\tau_j; y, u) + \int_{\tau_j}^T p_j(t) \Omega f(t, \tau_j; y, u) dt + \int_{\tau_j}^T p_j(t) [\frac{\partial g(t;...)}{\partial \tau_j} + \frac{\partial g(t;...)}{\partial y_j} D_j y(\tau_j^-)] dt +$$

$$+ \frac{\partial G(\boldsymbol{\tau}, \mathbf{y}, \mathbf{a})}{\partial \tau_j} + \frac{\partial G(\boldsymbol{\tau}, \mathbf{y}, \mathbf{a})}{\partial y_j} D_j y(\tau_j^-) + \sum_{\ell=j+1}^{N+1} \frac{\partial G(\boldsymbol{\tau}, \mathbf{y}, \mathbf{a})}{\partial y_\ell} [\tilde{\eta}_j(\tau_\ell^-) + \int_{\tau_j}^{\tau_\ell} R_j(\tau_\ell^-, t) \tilde{\eta}_j(t) dt]$$

--- (5.6)

The operator $\Omega$ applied to the Hamiltonian $h_j$ gives

$$\Omega h_j(\tau_j; y, u, p(.)) = h_j(\tau_j, y(\tau_j^-), a_{j-1}, p(\cdot)) - h_j(\tau_j, y(\tau_j^+), a_j, p(\cdot)) \equiv$$

$$\equiv \Omega F(\tau_j; y, u) + \int_{\tau_j}^T p_j(t) \Omega f(t, \tau_j; y, u) dt$$

--- (5.7)

thus

$$\frac{\partial J}{\partial \tau_j} = \Omega h_j(\tau_j; y, u, p(.)) + \int_{\tau_j}^{T} p_j(t)[\frac{\partial g(t,...)}{\partial \tau_j} + \frac{\partial g(t,...)}{\partial y_j} D_j y(\tau_j^-)] dt +$$

$$+ \frac{\partial G(\boldsymbol{\tau},\mathbf{y},\mathbf{a})}{\partial \tau_j} + \frac{\partial G(\boldsymbol{\tau},\mathbf{y},\mathbf{a})}{\partial y_j} D_j y(\tau_j^-) + \sum_{\ell=j+1}^{N+1} \frac{\partial G(\boldsymbol{\tau},\mathbf{y},\mathbf{a})}{\partial y_\ell}[\tilde{\eta}_j(\tau_\ell^-) + \int_{\tau_j}^{\tau_\ell} R_j(\tau_\ell^-,t)\tilde{\eta}_j(t)dt]$$

--- (5.8)

Our next task is to find Hamiltonian equations for the quantities $R_j(\tau_\ell^-, t)$ that appear in the last term on the right-hand side of (5.8). It follows from the results of section 3 that each $R_j(\tau_\ell^-, t)$ satisfies the following impulsive Volterra equation:

$$R_j(\tau_\ell^-, t) = \frac{\partial f(\tau_\ell, t, y(t), u(t))}{\partial y} + \sum_{k: t < \tau_k < T} \sum_{m=k+1}^{N} \frac{\partial g(\tau_\ell, \cdots)}{\partial y_m} \Gamma_{mk} \frac{\partial f(\tau_k, t, y(t), u(t))}{\partial y} +$$

$$+ \int_{t}^{\tau_\ell} R_j(\tau_\ell^-, s) \frac{\partial f(s, t, y(t), u(t))}{\partial y} ds +$$

$$+ \sum_{k: t < \tau_k < T} \sum_{m=k+1}^{N} \left( \int_{\tau_m}^{T} R_j(\tau_\ell^-, s) \frac{\partial g(s, \cdots)}{\partial y_m} ds \right) \Gamma_{mk} \frac{\partial f(\tau_k, t, y(t), u(t))}{\partial y}$$

--- (5.9)

Therefore, we define another set of Hamiltonians $H_\ell$, $1 \leq \ell \leq N+1$, in the following way:

$$H_\ell(t, y, u, r_\ell(.)) := f(\tau_\ell, t, y, u) + \int_{t}^{\tau_\ell} r_\ell(s) f(s, t, y, u) ds$$

--- (5.10)

Then the quantities $R_j(\tau_\ell^-, t)$ satisfy

$$R_j(\tau_\ell^-, t) = \frac{\partial H_\ell(t, y(t), u(t), R_j(\tau_\ell^-, \cdot))}{\partial y} +$$

$$+ \sum_{k: t < \tau_k < T} \sum_{m=k+1}^{N} \left( \frac{\partial g(\tau_\ell, \cdots)}{\partial y_m} + \int_{\tau_m}^{T} R_j(\tau_\ell^-, s) \frac{\partial g(s, \cdots)}{\partial y_m} ds \right) \Gamma_{mk} \frac{\partial f(\tau_k, t, y(t), u(t))}{\partial y},$$

for $\tau_j < t < \tau_\ell$

--- (5.11)

In the following theorem, we gather the results we have proved thus far in this section:

Theorem 5.1. The derivatives of the cost functional with respect to the impulsive times are given by (5.6), where the quantities $\eta_j$ are defined in (4.11), the co-states $p_j$ satisfy the Hamiltonian impulsive integral equations (5.5) with Hamiltonian functionals given by (5.4), and the quantities $R_j(\tau_\ell^-, t)$ satisfy the Hamiltonian impulsive integral equations (5.11) with Hamiltonian functionals given by (5.10). ///

By combining (5.8) with standard results of finite-dimensional optimization, we obtain:

Theorem 5.2. If there is an optimal collection of impulse times, say $\boldsymbol{\tau}^* = \{\tau_1^*, \tau_2^*, \cdots, \tau_N^*\}$ in the open set $S_N$, i.e. $0 < \tau_1^* < \tau_2^* < \cdots < \tau_N^* < T$, then the first order necessary condition for optimality is

$$\left( \Omega h_j(\tau_j; y, u, p(.)) + \int_{\tau_j}^{T} p_j(t) \left[ \frac{\partial g(t; \ldots)}{\partial \tau_j} + \frac{\partial g(t; \ldots)}{\partial y_j} D_j y(\tau_j^-) \right] dt + \right.$$

$$+ \frac{\partial G(\boldsymbol{\tau}, \mathbf{y}, \mathbf{a})}{\partial \tau_j} + \frac{\partial G(\boldsymbol{\tau}, \mathbf{y}, \mathbf{a})}{\partial y_j} D_j y(\tau_j^-) +$$

$$\left. + \sum_{l=j+1}^{N+1} \frac{\partial G(\boldsymbol{\tau}, \mathbf{y}, \mathbf{a})}{\partial y_l} [\eta_j(\tau_l^-) + \int_{\tau_j}^{\tau_l} R_j(\tau_l^-, t) \eta_j(t) dt] \right) \Bigg|_{\boldsymbol{\tau} = \boldsymbol{\tau}^*} = 0,$$

$\forall\, i = 1, 2, \ldots, N$

--- (5.12)

///

## 6. The case of impulsive ordinary differential equations.

When a controlled system is described by impulsive ordinary differential equations, the results of the previous sections still apply, but it is possible to carry out a number of simplifications. Because controlled impulsive ODEs are an important class of control systems, we examine the first-order necessary conditions for optimality with respect to variations of the impulse times. These results are, to the best of our knowledge, new for controlled impulsive ODEs, and they are not obtainable in any direct way from classical works on controlled impulsive ODEs.

Except where otherwise specified, we use the notation and conditions of section 2. We consider an impulsive controlled ODE system

$$\frac{dy(t)}{dt} = f(t, y(t), u(t)), \text{ for } t \notin \boldsymbol{\tau}; \quad y(0) = y_0;$$

$$y(\tau_i^+) = I(\tau_i, y(\tau_i^-), a_{i-1})$$

--- (6.1)

The problem (6.1) can be written in integral form as

$$y(t) = y_0 + \int_0^t f(s, y(s), u(s))\,ds + \sum_{i: 0 < \tau_i < t} I(\tau_i, y(\tau_i^-), a_{i-1})$$

--- (6.2)

Thus we have a particular case of (2.1) where the function f(t, s, y(s), u(s)) that appears in (2.1) is replaced by a function that does not depend on t in (6.1), and the function $g(t, \tau_{(t)}, y_{(t)}, a_{(t)})$ of (2.1) is replaced by a function that does not depend on t and is given by

$$g(\tau_{(t)}, y_{(t)}, a_{(t)}) = \sum_{i: 0 < \tau_i < t} I(\tau_i, y(\tau_i^-), a_{i-1})$$

--- (6.3)

In this case, we set

$$\psi_j(t) := \int_t^T p_j(s)\,ds$$

--- (6.4)

thus

$$\psi_j(T) = 0,\ p_j(t) = -\dot{\psi}_j(t),\ \int_{\tau_i}^{T} p_j(t)\,dt = \psi_j(\tau_i^+)$$

--- (6.5)

The Hamiltonians $h_j$ now become

$$h_j(t, y, u, \psi_j) = F(t, y, u) + \psi_j f(t, y, u)$$

--- (6.6)

Thus (5.5) becomes, in the context of the present section,

$$\dot{\psi}_j(t) =$$
$$= -\frac{\partial h_j(t, y(t), u(t))}{\partial y} - \sum_{\ell: t < \tau_\ell < T} \sum_{i=\ell+1}^{N} \psi_j(\tau_i^+) \frac{\partial I(\tau_\ell, y(\tau_\ell^-), a_{\ell-1})}{\partial y} \Gamma_{i\ell} \frac{\partial f(t, y(t), u(t))}{\partial y}$$

--- (6.7)

To calculate $\Gamma_{i\ell}$, we observe that the terms $\lambda_\ell(t)$, defined in (4.7), are now given by

$$\lambda_\ell(t) = \frac{\partial I(\tau_\ell, y(\tau_\ell^-), a_{\ell-1})}{\partial y} \equiv \lambda_\ell,\ \text{thus they are independent of t, and consequently}$$

$$\Gamma_{i\ell} = \sum_{\alpha} \sum_{\substack{(k_\alpha, \ldots, k_1): \\ i < k_\alpha < \cdots < k_1 < \ell}} \lambda_i \lambda_{k_\alpha} \cdots \lambda_{k_1}$$

--- (6.8)

For the calculation of the terms that correspond to $R_j(\tau_\ell^-, t)$, we set

$$\rho_{j\ell}(t) := \int_t^{\tau_\ell} R_j(\tau_\ell^-, s)\,ds,\ \text{for } t > \tau_j$$

--- (6.9)

and, by using reasoning analogous to that for the terms $\psi_j(t)$, we can express the Hamiltonians $H_\ell$ as

$$H_\ell(t,y,u,\rho) = f(t,y,u) + \rho f(t,y,u) \equiv H(t,y,u,\rho)$$

--- (6.10)

(thus each $H_\ell$ is actually independent of $\ell$), and the equations for $\rho_{j\ell}$ become

$$\dot{\rho}_{j\ell}(t) = -\frac{\partial H(t,y(t),u(t),\rho_{j\ell}(t))}{\partial y} -$$

$$-\sum_{k:t<\tau_k<T} \sum_{m=k+1}^{N} \left( \frac{\partial I(\tau_m, y(\tau_m^-), a_{m-1})}{\partial y} + \rho_{j\ell}(\tau_m^+)\frac{\partial I(\tau_m, y(\tau_m^-), a_{m-1})}{\partial y} \right) \cdot$$

$$\cdot \Gamma_{mk} \frac{\partial f(t,y(t),u(t))}{\partial y}, \quad \text{for } \tau_j < t < \tau_\ell$$

--- (6.11)

and it is seen that, apart from the condition $\tau_j < t < \tau_\ell$, the terms $\rho_{j\ell}(t)$ are in fact independent of j and $\ell$, and therefore we shall write $\rho(t)$ in lieu of $\rho_{j\ell}(t)$, and the range of values of t will be implied by the context of each particular equation. (According to the notational conventions stated in section 2, each $\rho(t)$ is an $n \times n$ matrix, whereas each $\eta_j(\cdot)$ or $\tilde{\eta}_j(\cdot)$ takes values in $R^n$.)

The terms $D_j y(\tau_j^-)$, defined in (4.4), become, in the case of impulsive ODEs,

$$D_j y(\tau_j^-) = f(\tau_j, y(\tau_j^-), a_{j-1}) + \frac{\partial I(\tau_j, y(\tau_j^-), a_{j-1})}{\partial \tau_j}$$

--- (6.12)

whereas the terms $\eta_j(t), \tilde{\eta}_j(t)$ become

$$\eta_j(t) \equiv \eta_j = f(\tau_j, y(\tau_j^-), a_{j-1}) - f(\tau_j, y(\tau_j^-) + I(\tau_j, y(\tau_j^-), a_{j-1}), a_j) +$$

$$+ \frac{\partial I(\tau_j, y(\tau_j^-), a_{j-1})}{\partial \tau_j} + \frac{\partial I(\tau_j, y(\tau_j^-), a_{j-1})}{\partial y} D_j y(\tau_j^-) ;$$

$$\tilde{\eta}_j(t) = \left(1 + \sum_{i:\ \tau_j < \tau_i < t} \sum_{\ell=j}^{i-1} \lambda_i \Gamma_{i\ell} \right) \eta_j$$

--- (6.13)

When all the information above is substituted into (5.8), we obtain

$$\frac{\partial J}{\partial \tau_j} = \Omega h_j(\tau_j; y, u, \psi_j) + \psi_j(\tau_j^+) \left[ \frac{\partial I(\tau_j, y(\tau_j^-), a_{j-1})}{\partial \tau_j} + \frac{\partial I(\tau_j, y(\tau_j^-), a_{j-1})}{\partial y} D_j y(\tau_j^-) \right] +$$

$$+ \frac{\partial G(\boldsymbol{\tau}, \mathbf{y}, \mathbf{a})}{\partial \tau_j} + \frac{\partial G(\boldsymbol{\tau}, \mathbf{y}, \mathbf{a})}{\partial y_j} D_j y(\tau_j^-) + \sum_{\ell=j+1}^{N+1} \frac{\partial G(\boldsymbol{\tau}, \mathbf{y}, \mathbf{a})}{\partial y_\ell} [\tilde{\eta}_j + \rho(\tau_j^+) \tilde{\eta}_j]$$

--- (6.14)

Therefore, the first-order necessary condition for a minimizer $\boldsymbol{\tau}^*$ in the open set $S_N$ is

$$\Omega h_j(\tau_j; y^*, u^*, \psi^*_j) + \psi^*_j(\tau_j^+) \left[ \frac{\partial I(\tau_j, y^*(\tau_j^-), a^*_{j-1})}{\partial \tau_j} + \frac{\partial I(\tau_j, y^*(\tau_j^-), a^*_{j-1})}{\partial y} D_j y^*(\tau_j^-) \right] +$$

$$+ \frac{\partial G(\boldsymbol{\tau}, \mathbf{y}^*, \mathbf{a}^*)}{\partial \tau_j} + \frac{\partial G(\boldsymbol{\tau}, \mathbf{y}^*, \mathbf{a}^*)}{\partial y_j} D_j y^*(\tau_j^-) + \sum_{\ell=j+1}^{N+1} \frac{\partial G(\boldsymbol{\tau}, \mathbf{y}^*, \mathbf{a}^*)}{\partial y_\ell} [\tilde{\eta}^*_j + \rho^*(\tau_j^+) \tilde{\eta}^*_j] \bigg|_{\boldsymbol{\tau}=\boldsymbol{\tau}^*} = 0$$

--- (6.15)

where the superscript * denotes optimality.

## 7. Outline of the necessary optimality conditions arising from variations in the values of the control.

A problem of optimal control for impulsive Volterra integral equations with fixed impulse times has been studied in [BS]. That paper contains (among other things) the first-order necessary conditions for optimality that arise from variations in the values of the control.

For the problem of the present paper, with variable impulse times, both the variation of the impulse times and the variation of the values of the control need to be analyzed. As far as the control action is concerned, the model used in the present paper is very similar, but not completely identical, to the model in [BS]. For this reason, we outline, in this section, without proofs, the first-order necessary conditions for optimality arising out of variations in the values of the control for the model of the present paper. The differences between [BS] and this section are such that the proofs would be practically the same as the proofs in [BS], although the final results are formally different and thus worth reporting.

The co-state (with respect to variations in the values of the control) is denoted by $\varphi(t)$ and it satisfies the following integral equation:

$$\varphi(t) = -\frac{\partial F(t, y(t), u(t))}{\partial y} - \sum_{i=1}^{N+1} \sum_{j=1}^{i} \frac{\partial G(\tau, y, a)}{\partial y_i} \Gamma_{ij} \frac{\partial f(\tau_j, t, y(t), u(t))}{\partial y} +$$

$$+ \int_t^T \varphi(s) \frac{\partial f(s, t, y(t), u(t))}{\partial y} ds +$$

$$+ \sum_{j: t < \tau_j < T} \sum_{i=j}^{N+1} \left( \int_t^T \varphi(s) \frac{\partial g(t, \tau_{(t)}, y_{(t)}, a_{(t)})}{\partial y_i} ds \right) \Gamma_{ij} \frac{\partial f(\tau_j, t, y(t), u(t))}{\partial y}$$

--- (7.1)

The derivatives of the cost functional with respect to $\mathbf{a}$ are given by

$$\frac{\partial J}{\partial a_i} = -\int_{\tau_i}^{\tau_{i+1}} \int_t^T \varphi(s) \frac{\partial f(s, t, y(t), a_i)}{\partial a_i} ds\, dt - \int_{\tau_i}^T \varphi(t) \frac{\partial g(t, \tau_{(t)}, y_{(t)}, a_{(t)})}{\partial a_i} dt +$$

$$+ \frac{\partial G(\tau, y, a)}{\partial a_i} + \int_{\tau_i}^{\tau_{i+1}} \frac{\partial F(t, y(t), a_i)}{\partial a_i} dt +$$

$$+ \sum_{k=i}^{N+1} \sum_{j=i}^{k} \frac{\partial G(\tau, y, a)}{\partial y_k} \Gamma_{kj} \left( \frac{\partial g(\tau_j, \tau_{(\tau_j)}, y_{(\tau_j)}, a_{(\tau_j)})}{\partial a_i} + \int_{\tau_i}^{\tau_{i+1}} \frac{\partial f(\tau_j, t, y(t), a_i)}{\partial a_i} dt \right)$$

--- (7.2)

Then the derivatives $\dfrac{\partial J}{\partial a_i}$ can be used in combination with standard optimization results to give the first order necessary conditions for optimality over all admissible values of **a**:

if $\mathbf{a}^*$ is optimal, then

$$\sum_{i=0}^{N}\left(\left.\frac{\partial J}{\partial a_i}\right|_{\mathbf{a}=\mathbf{a}^*}\right)\delta a_i \geq 0$$

--- (7.3)

for all admissible variations $\delta\mathbf{a} = \text{col}\{\delta a_i : 0 \leq i \leq N\}$ (i.e. the column vector with components $\delta a_i$) of the control $\mathbf{a}^*$. The inequality (7.3) is the variational extremum principle with respect to variations of the control.